\documentclass[11pt]{article}
\usepackage{amsfonts,amsmath,amssymb,amssymb,tikz,color}
\usetikzlibrary{arrows}
\usetikzlibrary{decorations.pathreplacing}

\newtheorem{theorem}{\bf Theorem}[section]
\newtheorem{corollary}[theorem]{\bf Corollary}
\newtheorem{lemma}[theorem]{\bf Lemma}

\newtheorem{remark}[theorem]{\bf Remark}
\newtheorem{problem}[theorem]{\bf Problem}

\newcommand{\proof}{\noindent{\bf Proof.\ }}
\newcommand{\qed}{\hfill $\square$ \bigskip}

\newcommand{\Tr}{{\rm Tr}}
\newcommand{\ecc}{{\rm ecc}}

\begin{document}

\title{\bf New  transmission irregular chemical graphs}

\author{
Kexiang Xu $^{a,b}$
\and
Jing Tian $^{a,b}$
\and
Sandi Klav\v zar  $^{c,d,e}$
}

\date{}
\maketitle

\begin{center}
$^a$ College of Mathematics, Nanjing University of  Aeronautics \& Astronautics\\ Nanjing, Jiangsu, 210016, PR China\\
 \texttt{kexxu1221@126.com} \\
 \texttt{jingtian526@126.com} \\
\medskip

$^b$ MIIT Key Laboratory of Mathematical Modelling and High Performance\\
 Computing of Air Vehicles, Nanjing, Jiangsu, 210016, China\\
\medskip

$^c$ Faculty of Mathematics and Physics, University of Ljubljana, Slovenia\\
\texttt{sandi.klavzar@fmf.uni-lj.si}\\
\medskip

$^d$ Institute of Mathematics, Physics and Mechanics, Ljubljana, Slovenia\\
\medskip

$^e$ Faculty of Natural Sciences and Mathematics, University of Maribor, Slovenia\\
\medskip

\end{center}

\begin{abstract}
 The transmission of a vertex $v$ of a  (chemical) graph $G$ is the sum of distances from $v$ to  other vertices in $G$. If any two vertices of $G$ have different transmissions, then $G$ is a transmission irregular graph. It is shown that for any odd number $n\geq 7$ there exists a transmission irregular chemical tree of order $n$. A construction is provided which generates new transmission irregular (chemical) trees. Two additional families of chemical graphs are characterized by property of transmission irregularity and two sufficient condition provided which guarantee that the transmission irregularity is preserved upon adding a new edge.
\end{abstract}


\medskip\noindent
{\bf AMS Subj. Class:} 05C05, 05C12, 05C92

\section{Introduction}
\label{sec:intro}

In chemical graph theory, molecules are naturally represented by (chemical) graphs. In the next step, the graph distance function is an obvious tool for exploring chemical graphs, which in turn reflect the physico-chemical properties of the corresponding (organic) compounds, cf.~\cite{RoKi02}.  The Wiener index~\cite{wiener-1947} is a famous example, but many other distance-based (possibly combined with vertex degrees) topological indices have been studied such as Schultz index~\cite{Sc89, liu-2022},
 hyper-Wiener index~\cite{Ra93},
 Gutman index~\cite{Gu94, liu-2022},
 vertex-Szeged index~\cite{Gu94},
 PI index~\cite{KhKaAg01},
 edge-Wiener index~\cite{KhYoAsWa09},
 edge-Szeged index~\cite{GuAs08},
 Wiener polarity index~\cite{ali-2021, haoer-2021},
and more.
The area is still very active; for a survey on graphs extremal with respect to distance-based topological indices see~\cite{XLDGF2014}, and for a selection of recent developments with a focus on applications see~\cite{arock-2022, arock-2021, brezo-2021, das-2021, shira-2019, sorgun-2022, yu-2022, zhang-2022}.

Exploring all these indices can be interesting from a mathematical point of view, but it is also  much important from a chemical point of view, as it turns out in practice that several indices need to be combined to determine the properties of molecules. Moreover, this approach has found applications in many other areas such as  communication theory, facility location, crystallography, and even in ornithology. As a point of interest for the latter we mention that it was shown in~\cite{BaCa08} that  the interaction between a flock of birds depends more intimately on the topological distance rather than the Euclidian distance.

 If $G = (V(G), E(G))$ is a graph and $x,y\in V(G)$, then $d_G(x,y)$ denotes the shortest-path distance between $x$ and $y$  in $G$. The sum of all distances from a vertex $x$ to other vertices is a basic building block in exploration of metric properties of a graph and is called the {\em transmission} of $x$ and denoted by $\Tr_G(x)$. That is,
$$\Tr_G(x)=\sum\limits_{u\in V(G)\setminus\{x\}}d_G(x,u)\,.$$
 That the transmission of a vertex is indeed a fundamental concept is demonstrated by the fact that it is also known by several other names, such as the status of a vertex~\cite{abiad-2021, QZ2020} and the total distance of a vertex~\cite{cava-2019, kla-2013}.

 The \textit{transmission set} of $G$ is
$$\Tr(G) = \{\Tr_G(x):\ x\in V(G)\}\,.$$
If $|\Tr(G)|=n(G)$ holds, where $n(G)$ denotes the order of $G$, then $G$ is \textit{transmission irregular}, TI for short. Recalling that the Wiener complexity of a $G$ is the number of different transmission of its vertices~\cite{AAKS2014}, see also~\cite{dob-2022, xu-2021}, we can equivalently say that TI graphs are the graphs with maximum Wiener complexity.

 Since almost no graph is transmission irregular~\cite{AK2018}, the search for such graphs has become of interest to several groups of researchers. Results to date have been presented in~\cite{al-yakoob-2020, al-yakoob-2022a, al-yakoob-2022b, bezh-2021, bezh-2022, dobrynin-2019, dobrynin-2019b, dobrynin-2019c, DS2020, XK21}. In this paper we continue this line of research with a focus on chemical graphs. In the next section we list definitions needed, recall some known results, and prove a couple of results to be used later. In Section~\ref{sec:tree} we investigate transmission irregularity of chemical trees while in Section~\ref{sec:cycle} we consider families of chemical graphs containing a few short cycles. We conclude the paper with some open problems.

\section{Preliminaries}
\label{sec:prelim}

All graphs considered in this paper are finite, simple and, unless stated otherwise, also connected.  For $X\subseteq V(G)$, let $G-X$ be the subgraph of $G$ obtained from $G$ by removing the vertices from $X$ {and the edges incident with them}, in particular, $G - \{v\}$ will be briefly denoted by $G-v$. Similarly, for $F\subseteq E(G)$, $G-F$ is the spanning subgraph of $G$ obtained by removing the edges of $F$ and if $e\in E(G)$ then we will write $G-e$ for $G-\{e\}$. 
The \textit{eccentricity} $\ecc_G(v)$ of a vertex $v\in V(G)$ is the maximum distance from $v$ to all other vertices in $G$. If $uv\in E(G)$, then $n_u$ (or $n_G(u)$ if the graph $G$ is necessarily mentioned) is the number of vertices in $G$ closer to $u$ than to $v$ and $n_v$ (or $n_G(v)$ for completeness) is similarly defined.

 A vertex $v$ with $\deg_G(v)=1$ is called a \textit{pendant vertex} (also {\em leaf} when $G$ is a tree) in $G$, and the edge incident with a pendant vertex is
called a \textit{pendant edge}. A path $P:=u_ku_{k-1}\cdots u_2u_1$ with natural adjacency relation in a graph $G$ is a \textit{proper pendant path} in $G$ if $\deg_G(u_k)\geq 3$, $\deg_G(u_1)=1$, and $\deg_G(u_i)=2$ for $i\in \{2,3,\ldots,k-1\}$, where $u_{k}$ is its root. If both $u_k$ and $u_1$ in $P$ have degrees at least $3$ and each of $u_j$ with $j\in \{2,3,\ldots,k-1\}$ has degree $2$, then $P$ is an \textit{internal path} in $G$ with two \textit{terminals} $u_k$ and $u_1$. Specially, if $u_1$ and $u_k$ have degrees at least $2$, then the above $P$ is a \textit{weak internal path} with two weak terminals $u_1$ and $u_k$.

 The definition of a chemical graph is still a matter of debate, but we will stick to the most common and simple one: A graph $G$ is a {\em chemical graph} if its maximum degree is at most $4$. A vertex of degree at least $3$ is a \textit{branching vertex}. A tree with a unique branching vertex $v$ is \textit{starlike}. A starlike tree $T$ with branching vertex $v$ will be denoted by $T=T(n_1,\ldots,n_k)$ if $T-v$ consists of $k$ disjoint paths of orders $n_1,\ldots, n_k$, respectively. And the pendant path of length $n_i$ from $v$ is called an {\em $n_i$-arm} in it.

 For an induced subgraph $H$ of $G$, we say that the transmission set of $H$ in $G$ is $\Tr_G(H)=\{\Tr_G(u):u\in V(H)\}$. In particular, $\Tr_G(G )= \Tr(G)$. For an induced subgraph $H$ of a graph $G$, if $|\Tr_G(H)|=n(H)$, then $H$ is a \textit{partially transmission irregular} subgraph of $G$.

For a positive integer $k$ we use the notation $[k]=\{1,\ldots, k\}$ and $[k]_0=\{0,1,\ldots, k\}$. For a set $A$ of integers and $i\in \mathbb{Z}$, we denote by $A+i$ the usual coset, that is, $A+i=\{a+i:a\in A\}$.

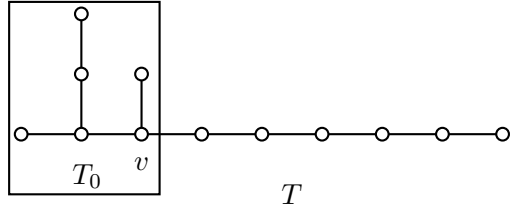
\begin{figure}[ht!]
\begin{center}
\begin{tikzpicture}[scale=0.8,style=thick]
\def\vr{3pt} 

\draw (0,2) -- (0,1) -- (1,1) -- (2,1) -- (3,1) -- (4,1) -- (5,1) -- (6,1);
\draw (-1,3) -- (-1,2) -- (-1,1) -- (0,1);
\draw (-2,1) -- (-1,1);

\draw (1,1)  [fill=white] circle (\vr);
\draw (0,1)  [fill=white] circle (\vr);
\draw (2,1)  [fill=white] circle (\vr);
\draw (3,1)  [fill=white] circle (\vr);
\draw (4,1)  [fill=white] circle (\vr);
\draw (5,1)  [fill=white] circle (\vr);
\draw (6,1)  [fill=white] circle (\vr);
\draw (0,2)  [fill=white] circle (\vr);

\draw (-1,1)  [fill=white] circle (\vr);
\draw (-1,2)  [fill=white] circle (\vr);
\draw (-1,3)  [fill=white] circle (\vr);
\draw (-2,1)  [fill=white] circle (\vr);

\draw (-2.2,0) rectangle (0.3,3.2);


\draw (2.5,0) node {$T$};

\draw (0,0.5) node {$v$};

\draw (-0.9,0.3) node {$T_0$};

\end{tikzpicture}
\end{center}
\caption{Tree $T$ with a DBTM  subtree $T_0$ rooted at $v$.}
\label{fig: DBTM}
\end{figure}

For any tree $T$ and its subtree $T_0$ with a non-leaf vertex $v\in V(T_0)$, we denote by $V_j$ the set of vertices at distance $j$ from $v$ in $T_0$. Let $a = \ecc_{T_0}(v)$. Then $V(T^*)=\cup_{j=1}^{a}V_j$ is a \textit{distance-based partition} of the forest $T_0-v$ at $v$. If $\min\limits_{u\in V_{j+1}}\Tr_T(u)\geq \max\limits_{u\in V_j}\Tr_T(u)$ for any $j\in [a-1]$ in the above partition, then $T_0$ is a \textit{distance-based transmission monotonic (DBTM for short) subtree} of $T$ at $v$. See an example of a DBTM subtree in Fig.~\ref{fig: DBTM}. In particular, if $T_0=T$, then $T$ is a DBTM tree. If $\min\limits_{u\in V_{j+1}}\Tr_T(u)\geq \max\limits_{u\in V_{j-1}}\Tr_T(u)$ for any $j\in [a-1]\setminus [1]$ in the above distance-based partition of $T$ at $v\in V(T)$, then $T$ is a $2$-DBTM tree at $v$. If $T$ is a DBTM tree at $v\in V(T)$, then $\min\limits_{u\in V_{j+1}}\Tr_T(u)\geq \max\limits_{u\in V_j}\Tr_T(u)\geq \min\limits_{u\in V_j}\Tr_T(u)\geq \max\limits_{u\in V_{j-1}}\Tr_T(u)$, which implies that $T$ is a $2$-DBTM tree at $v$. Therefore DBTM tree is a special $2$-DBTM tree with the same root.

A set of positive integers is odd (even, resp.) if it consists of odd (even, resp.) integers. A family $A = \cup_{i=1}^{k}A_i$ of sets of positive integers has \textit{intersecting parity} if $A_p$ and $A_{p+1}$ have different parities for any $p\in [k-1]$. Moreover, similarly as above, if $\min A_{j+1}\geq \max A_{j-1}$ for any $j\in [k-1]\setminus [1]$, then $A$ is \textit{$2$-distance monotonic}.

\begin{lemma}\label{2dm}
Let $A = \cup_{i=1}^{k}A_i$ be a $2$-distance monotonic family of sets of positive integers. If $A$ has intersecting parity, then  the sets $A_i$ are pairwise disjoint.
\end{lemma}

\proof Without loss of generality, we assume that $A_1$ is odd. From the condition that $A$ has intersecting parity, we deduce that $A_j$ has the same parity with its subscript $j$ for any $j\in [k]$. Then $A_p\cap A_q=\emptyset$ if  $p,q\in [k]$ have different parities. It follows that $A^{(1)}\cap A^{(2)}=\emptyset$ where $A^{(1)}$ is the union of sets $A_i$ with odd $i\in [k]$  and $A^{(2)}$ is the union of sets $A_i$ with even $i\in [k]$. In view of the $2$-distance monotonic property of $A$, we conclude that both $A^{(1)}$ and $A^{(2)}$ are pairwise disjoint. Thus the result follows immediately. \qed

From Lemma~\ref{2dm}, the following result is obvious.

\begin{corollary}\label{co-2dm}
Let $A = \cup_{i=1}^{k}A_i$ be a $2$-distance monotonic family of sets of positive integers, and let $t$ be an even positive integer. If $A$ has intersecting parity, then both $\cup_{i=1}^{k}(A_i+it)$ and $\cup_{i=1}^{k}(A_i+a)$ are pairwise disjoint, where $a$ is a constant. \end{corollary}

\begin{lemma}\label{equal} {\rm \cite{Bala}}
If $G$  is a graph with $n(G) > 2$ and $uv\in E(G)$, then $\Tr(u)-\Tr(v)=n_v-n_u$.
\end{lemma}

\begin{lemma}\label{pend} {\rm \cite{XK21}}
Let $G$ be a graph with  $n(G)=n$ and $v$ a vertex with $\deg(v)\geq 3$. If $P=uv_1v_2\cdots v_{x-1}v$ is a pendant path with natural adjacency relation attaching at $v$, where $\deg(u)=1$ and $x<\frac{n}{2}$, then $\Tr(v_{x-1}) - \Tr(v)=n-2x$.
\end{lemma}

\begin{lemma}\label{internal}
Let $G$ be a graph with  $n(G)=n$ and $P=vv_1v_2\cdots v_kv^*$ is a weak internal path in $G$ with two weak terminals $v$ and $v^*$ such that each edge in $P$ is a cut edge. If $\Tr(v_1)-\Tr(v)=a>0$, then $\Tr(v_{j})-\Tr(v)=j(a+j-1)$ for any $j\in [k]$.
\end{lemma}
\proof By Lemma \ref{equal}, we have $n_{v}-n_{v_1}=a$. Since each edge in $P$ is a cut edge, we have $n_{v_1}-n_{v_2}=a+2$, $n_{v_2}-n_{v_3}=a+4$, $\ldots$, $n_{v_{j-1}}-n_{v_j}=a+2(j-1)$. From Lemma \ref{equal}, we have $\Tr(v_2)-\Tr(v_1)=a+2$, $\Tr(v_3)-\Tr(v_2)=a+4$, $\ldots$, $\Tr(v_{j})-\Tr(v_{j-1})=a+2(j-1)$. Note that $v_{p-1}=v$ if $p=1$.  It follows that \begin{align*}
  \Tr(v_j)-\Tr(v) & = \sum\limits_{p=1}^{j}\Big(\Tr(v_p)-\Tr(v_{p-1})\Big) \\
   & = j(a+j-1),
\end{align*} completing the proof.
\qed

\section{ Transmission irregular chemical trees}
\label{sec:tree}

\begin{figure}[ht!]
\begin{center}
\begin{tikzpicture}[scale=1.0,style=thick]
\tikzstyle{every node}=[draw=none,fill=none]
\def\vr{3pt} 

\begin{scope}[yshift = 0cm, xshift = 0cm]
\path (0,0) coordinate (x1);
\path (1,0) coordinate (x2);
\path (3,0) coordinate (x3);
\path (0.7,0.7) coordinate (x4);
\path (2,2) coordinate (x5);
\path (0.7,-0.7) coordinate (x6);
\path (2,-2) coordinate (x7);
\path (3.7,0.7) coordinate (x8);
\path (5,2) coordinate (x9);
\path (3.7,-0.7) coordinate (x10);
\path (5,-2) coordinate (x11);
\draw (x1) -- (x2);
\draw (x1) -- (x4);
\draw (x1) -- (x6);
\draw (x3) -- (x8);
\draw (x3) -- (x10);
\draw (1.0,0) -- (1.5,0);
\draw (2.5,0) -- (3.0,0);

\draw (0.7,0.7) -- (1.0,1.0);
\draw (1.7,1.7) -- (2,2);
\draw (3.7,0.7) -- (4.0,1.0);
\draw (4.7,1.7) -- (5,2);

\draw (0.7,-0.7) -- (1.0,-1.0);
\draw (1.7,-1.7) -- (2,-2);
\draw (3.7,-0.7) -- (4.0,-1.0);
\draw (4.7,-1.7) -- (5,-2);

\draw (x1)  [fill=white] circle (\vr);
\draw (x2)  [fill=white] circle (\vr);
\draw (x3)  [fill=white] circle (\vr);
\draw (x4)  [fill=white] circle (\vr);
\draw (x5)  [fill=white] circle (\vr);
\draw (x6)  [fill=white] circle (\vr);
\draw (x7)  [fill=white] circle (\vr);
\draw (x8)  [fill=white] circle (\vr);
\draw (x9)  [fill=white] circle (\vr);
\draw (x10)  [fill=white] circle (\vr);
\draw (x11)  [fill=white] circle (\vr);
\draw (2.1,0) node {$\cdots$};
\draw (1.2,1.2) node {$\cdot$};
\draw (1.35,1.35) node {$\cdot$};
\draw (1.5,1.5) node {$\cdot$};
\draw (4.2,1.2) node {$\cdot$};
\draw (4.35,1.35) node {$\cdot$};
\draw (4.5,1.5) node {$\cdot$};
\draw (1.2,-1.2) node {$\cdot$};
\draw (1.35,-1.35) node {$\cdot$};
\draw (1.5,-1.5) node {$\cdot$};
\draw (4.2,-1.2) node {$\cdot$};
\draw (4.35,-1.35) node {$\cdot$};
\draw (4.5,-1.5) node {$\cdot$};
\draw (2.05,0.6) node {$k$};
\draw (1.0,1.8) node {$a_1$};
\draw (4.0,1.8) node {$b_1$};
\draw (1.9,-1.1) node {$a_2$};
\draw (4.9,-1.1) node {$b_2$};
\draw [decorate, decoration = {brace}] (0.9,0.2) --  (3.1,0.2);
\draw [decorate, decoration = {brace}] (0.5,0.8) --  (2.0,2.3);
\draw [decorate, decoration = {brace}] (3.5,0.8) --  (5.0,2.3);
\draw [decorate, decoration = {brace}] (0.8,-0.5) --  (2.3,-2);
\draw [decorate, decoration = {brace}] (3.8,-0.5) --  (5.3,-2);
\end{scope}
\end{tikzpicture}
\end{center}
\caption{The tree $H^{k}(a_1,a_2;b_1,b_2)$}
\label{fig:H-k}
\end{figure}
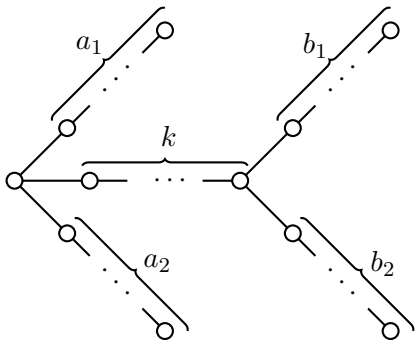

In~\cite{al-yakoob-2020, XK21} some TI starlike trees are determined, in particular, the TI starlike trees with maximum degree $3$ are characterized in~\cite{al-yakoob-2020} with a complicated condition. It is proved in~\cite{XK21} that $T=T(a,a+1,\ldots,a+k)$ is TI if $n(T)$ is odd. A tree $H^{k}(a_1,a_2;b_1,b_2)$ is obtained by attaching two pendant paths of lengths $b_1$, $b_2$, respectively, at a leaf on the $k$-arm of $T(a_1,a_2,k)$,
 see Fig.~\ref{fig:H-k}.

Before stating our  first main result, we make the following comment which will be used frequently in the subsequent proof.

\begin{remark}\label{Re1} Let $a<b$ be two positive integers. If $b^2-a^2$ is even (odd, resp.), then $b-a$ and  $b+a$ are both even (odd, resp.). \end{remark}

\proof Clearly, $b^2-a^2=(b+a)(b-a)$. Note that $b+a$ and $b-a$ have a same parity since $(b+a)+(b-a)=2b$ is even. Thus our result follows immediately.
\qed

By a computer search we find that there is no TI tree of order at most $6$. On the other hand, we have the following result.

\begin{theorem}\label{C-TI}
If $n\geq 7$ is an odd integer, then there exists a  TI chemical tree of order $n$.
\end{theorem}

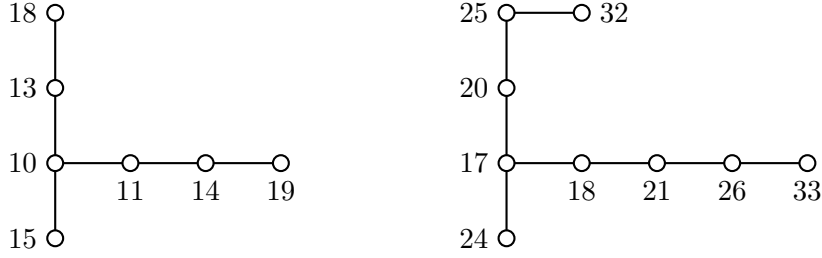
\begin{figure}[ht!]
\begin{center}
\begin{tikzpicture}[scale=1.0,style=thick]
\tikzstyle{every node}=[draw=none,fill=none]
\def\vr{3pt} 

\begin{scope}[yshift = 0cm, xshift = 0cm]
\path (0,0) coordinate (v1);
\path (1,0) coordinate (v2);
\path (2,0) coordinate (v3);
\path (3,0) coordinate (v4);
\path (0,-1) coordinate (v5);
\path (0,1) coordinate (v6);
\path (0,2) coordinate (v7);
\draw (v1) -- (v2) -- (v3) -- (v4);
\draw (v5) -- (v1) -- (v6) -- (v7);
\draw (v1)  [fill=white] circle (\vr);
\draw (v2)  [fill=white] circle (\vr);
\draw (v3)  [fill=white] circle (\vr);
\draw (v4)  [fill=white] circle (\vr);
\draw (v5)  [fill=white] circle (\vr);
\draw (v6)  [fill=white] circle (\vr);
\draw (v7)  [fill=white] circle (\vr);
\draw[left] (v5)++(-0.1,0.0) node {$15$};
\draw[left] (v1)++(-0.1,0.0) node {$10$};
\draw[left] (v6)++(-0.1,0.0) node {$13$};
\draw[left] (v7)++(-0.1,0.0) node {$18$};
\draw[below] (v2)++(0.0,-0.1) node {$11$};
\draw[below] (v3)++(0.0,-0.1) node {$14$};
\draw[below] (v4)++(0.0,-0.1) node {$19$};
\end{scope}

\begin{scope}[yshift = 0cm, xshift = 6cm]
\path (0,0) coordinate (v1);
\path (1,0) coordinate (v2);
\path (2,0) coordinate (v3);
\path (3,0) coordinate (v4);
\path (0,-1) coordinate (v5);
\path (0,1) coordinate (v6);
\path (0,2) coordinate (v7);
\path (1,2) coordinate (v9);
\path (4,0) coordinate (v8);
\draw (v1) -- (v2) -- (v3) -- (v4);
\draw (v5) -- (v1) -- (v6) -- (v7);
\draw (v7) -- (v9);
\draw (v4) -- (v8);
\draw (v1)  [fill=white] circle (\vr);
\draw (v2)  [fill=white] circle (\vr);
\draw (v3)  [fill=white] circle (\vr);
\draw (v4)  [fill=white] circle (\vr);
\draw (v5)  [fill=white] circle (\vr);
\draw (v6)  [fill=white] circle (\vr);
\draw (v7)  [fill=white] circle (\vr);
\draw (v8)  [fill=white] circle (\vr);
\draw (v9)  [fill=white] circle (\vr);
\draw[left] (v5)++(-0.1,0.0) node {$24$};
\draw[left] (v1)++(-0.1,0.0) node {$17$};
\draw[left] (v6)++(-0.1,0.0) node {$20$};
\draw[left] (v7)++(-0.1,0.0) node {$25$};
\draw[below] (v2)++(0.0,-0.1) node {$18$};
\draw[below] (v3)++(0.0,-0.1) node {$21$};
\draw[below] (v4)++(0.0,-0.1) node {$26$};
\draw[right] (v9)++(0.1,0.0) node {$32$};
\draw[below] (v8)++(0.0,-0.1) node {$33$};
\end{scope}

\end{tikzpicture}
\end{center}
\caption{Sporadic TI chemical trees}
\label{fig:7-and-9}
\end{figure}

\proof
 TI chemical trees of order $7$ and $9$ are displayed in Fig.~\ref{fig:7-and-9}, where, for each vertex, we also give its transmission.

In the rest we assume that $n\geq 11$ is odd.

If $n=4a+3$, then $a\geq 2$. Now we consider the tree $T=H^2(a-1,a;a,a+1)$ and prove that $T$ is TI. Let $u$ be the vertex of degree $3$ in $T$ at which two pendant paths of lengths $a$, $a+1$, respectively,  are attached. Assume that $\Tr_T(u)=x$ and $uv\in E(T)$ with $\deg_T(v)=2$ and $vw\in E(T)$. Then $w$ is the other vertex of degree $3$ in $T$ with $\Tr_T(v)=x+1$, and $\Tr_T(w)=x+4$.  Let $A_1$ and $A_3$ be the sets of transmissions of vertices on the pendant paths of lengths $a-1$ and $a$, respectively, attached at $w$,  and let $A_2$ and $A_4$ be the sets of transmissions of vertices on the pendant paths of lengths $a$ and $a+1$, respectively, attached at $u$. From the structure of $T$ and Lemma~\ref{internal} we have
\begin{align*}
A_1 & = \{2ka+(k+2)^2:\ k\in[a-1]\}+x, \\
A_2 & = \{2ka+(k+1)^2-1:\ k\in[a]\}+x, \\
A_3 & = \{2ka+(k+1)^2-1:\ k\in[a]\}+4+x, \\
A_4 & = \{2ka+k^2:\ k\in[a+1]\}+x.
\end{align*}
 Next we prove that $A_i\cap A_j=\emptyset$ for any $i,j\in [4]$. If $2k_1a+(k_1+2)^2=2k_2a+(k_2+1)^2-1$ with $k_1\in [a-1]$ and $k_2\in [a]$, then $k_1<k_2$. It follows that $2(k_2-k_1)a-1=(k_1+2)^2-(k_2+1)^2$. By Remark \ref{Re1}, $k_1-k_2+1$ is odd, that is, $k_1-k_2\leq -2$ is even. However, we have $2(k_2-k_1)a-1>0>(k_1+k_2+3)(k_1-k_2+1)$ as a contradiction. Thus $A_1\cap A_2=\emptyset$. Similarly as above, we have $A_1\cap A_3=\emptyset$. If $2k_1a+(k_1+2)^2=2k_2a+k_2^2$ with $k_1\in [a-1]$ and $k_2\in [a+1]$, then $2(k_2-k_1)a=(k_1+2)^2-k_2^2$ with $k_2>k_1$. Note that $k_1-k_2\leq -2$ is even from Remark~\ref{Re1}. But $2(k_2-k_1)a>0\geq(k_1+k_2+2)(k_1-k_2+2)$ is a clear contradiction. Therefore $A_1\cap A_4=\emptyset$.  Note that $A_3=A_2+4$. Then $A_2\cap A_3=\emptyset$ since $|s-t|>2a\geq 4$ for any $s,t\in A_2$. If $2k_1a+(k_1+1)^2-1=2k_2a+k_2^2$ with $k_1\in [a]$ and $k_2\in [a+1]$, then $(k_1+1)^2-k_2^2=2(k_2-k_1)a+1$ with $k_2>k_1$, which implies that $k_1-k_2\leq -2$ is even from Remark \ref{Re1}. But $2(k_2-k_1)a+1>0>(k_1+k_2+1)(k_1-k_2+1)$ occurs contradictorily. Thus $A_2\cap A_4=\emptyset$. Similarly as above, we have $A_3\cap A_4=\emptyset$ as desired.

If $n=4b+1$, then $b\geq 3$. We consider the tree $T=H^2(b-2,b-1;b,b+1)$ and prove that $T$ is TI. Let $z_1$ be the vertex of degree $3$ in $T$ at which two pendant paths of lengths $b$, $b+1$, respectively, are attached, $z_1z_2\in E(T)$ with $\deg_T(z_2)=2$ and $z_2z_3\in E(T)$. Assume that $\Tr_T(z_1)=y$. Then, by Lemma \ref{equal}, we have $\Tr_T(z_2)=y+3$ and $\Tr_T(z_3)=y+8$. Denote by $B_i$ the set of transmissions of vertices on the pendant path of length $b+2-i$ attached at $z_1$ or $z_3$ with $i\in [4]$. From the structure of $T$ and Lemma~\ref{internal} we have
\begin{align*}
B_1 & = \{2ka+k^2-2k:\ k\in[b+1]\}+y, \\
B_2 & = \{2ka+k^2:\ k\in[b]\}+y, \\
B_3 & = \{2ka+(k+1)^2-1:\ k\in[b-1]\}+8+y, \\
B_4 & = \{2ka+(k+2)^2-4:\ k\in[b-2]\}+8+y.
\end{align*}
Now it suffices to prove that $B_i\cap B_j=\emptyset$ for any $i,j\in [4]$. If $2k_1b+k_1^2-2k_1=2k_2b+k_2^2$ with $k_1\in [b+1]$ and $k_2\in [b]$, then $2(k_1-k_2)b-1=k_2^2-(k_1-1)^2$ with $k_1>k_2$. By Remark~\ref{Re1}, $k_2-k_1\leq -2$ is even, which implies that $2(k_1-k_2)b-1>0>(k_2-k_1+1)(k_2+k_1-1)$ as a contradiction. This yields that $B_1\cap B_2=\emptyset$. If $2k_1b+k_1^2-2k_1=2k_2b+(k_2+1)^2+7$ with $k_1\in [b+1]$ and $k_2\in [b-1]$, then $2(k_1-k_2)b-8=(k_2+1)^2-(k_1-1)^2$ with $k_1>k_2$. Using Remark~\ref{Re1} again, we find that $k_2-k_1\leq -2$ is even. If $k_2-k_1=-2$, then we have $b=2$, contradicting to the fact that $b\geq 3$. If $k_2-k_1\leq -4$, we have $2(k_1-k_2)b-8\geq 8b-8>0>(k_2-k_1+2)(k_2+k_1)$ as a contradiction, again. Thus we get $B_1\cap B_3=\emptyset$. If $2k_1b+k_1^2-2k_1=2k_2b+(k_2+2)^2+4$ with $k_1\in [b+1]$ and $k_2\in [b-2]$, then $2(k_1-k_2)b-5=(k_2+2)^2-(k_1-1)^2$ with $k_1>k_2$. By Remark~\ref{Re1}, we observe that $k_2-k_1\leq -2$ is even. If $k_2-k_1=-2$, then $2(k_1-k_2)b-5=4b-5=2k_2+3$, that is, $k_2=4b-8>b-2$, contradicting to the fact $k_2\in [b-2]$ with $b\geq 3$. If $k_2-k_1\leq-4$, then $2(k_1-k_2)b-5\geq 4b-5>0>(k_2+2)^2-(k_1-1)^2$ as a clear contradiction again. Therefore $B_1\cap B_4=\emptyset$. If $2k_1b+k_1^2=2k_2b+(k_2+1)^2+7$ with $k_1\in [b]$ and $k_2\in [b-1]$, then $2(k_1-k_2)b-7=(k_2+1)^2-k_1^2$ with $k_1>k_2$, which implies that $k_2-k_1\leq -2$ is even from Remark~\ref{Re1}. But we deduce that $2(k_1-k_2)b-7\geq 4b-7>0>(k_2+1)^2-k_1^2$ as a contradiction. Thus $B_2\cap B_3=\emptyset$. Similarly as above, we have $B_2\cap B_4=\emptyset$. If $2k_1b+(k_1+1)^2+7=2k_2b+(k_2+2)^2+4$ with $k_1\in [b-1]$ and $k_2\in [b-2]$, then we have $2(k_1-k_2)b+3=(k_2+2)^2-(k_1+1)^2$ with $k_1>k_2$. Taking Remark~\ref{Re1} into account, we observe that $k_2-k_1\leq -2$ is even. However, $2(k_1-k_2)b+3\geq 4b+3>0>(k_2+2)^2-(k_1+1)^2$ is a contradiction, again, implying that $B_3\cap B_4=\emptyset$.
\qed

We next provide a method for constructing a TI tree from a tree with a DBTM subtree.

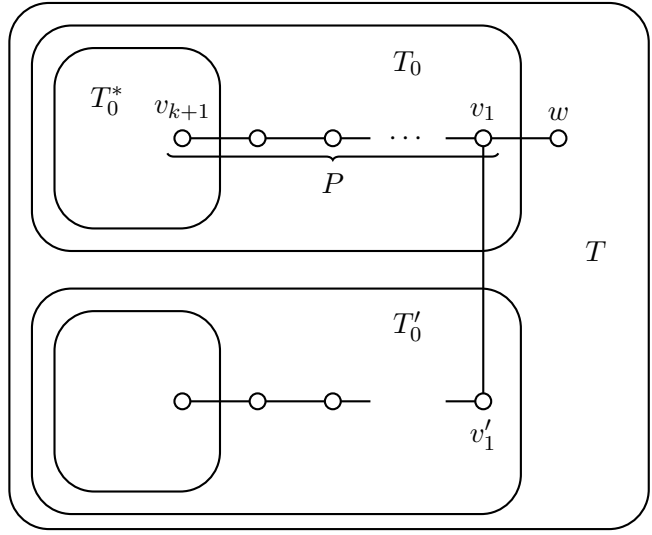
\begin{figure}[ht!]
\begin{center}
\begin{tikzpicture}[scale=1.0,style=thick]
\tikzstyle{every node}=[draw=none,fill=none]
\def\vr{3pt} 

\begin{scope}[yshift = 0cm, xshift = 0cm]
\path (6,1.5) coordinate (v1);
\path (4,1.5) coordinate (vk-1);
\path (3,1.5) coordinate (vk);
\path (2,1.5) coordinate (vk+1);
\path (7,1.5) coordinate (w);
\draw (vk+1) -- (vk) -- (vk-1);
\draw (vk-1) -- (4.5,1.5);
\draw (5.5,1.5) -- (v1);
\draw (v1) -- (6, -2);
\draw (v1) -- (w);
\draw (v1)  [fill=white] circle (\vr);
\draw (vk-1)  [fill=white] circle (\vr);
\draw (vk)  [fill=white] circle (\vr);
\draw (vk+1)  [fill=white] circle (\vr);
\draw (w)  [fill=white] circle (\vr);
\draw[rounded corners=15pt] (0,0) rectangle ++(6.5,3);
\draw[rounded corners=15pt] (0.3,0.3) rectangle ++(2.2,2.4);
\draw[rounded corners=15pt] (-0.3,-3.7) rectangle ++(8.5,7.0);
\draw (5.0,1.5) node {$\cdots$};
\draw[above] (vk+1)++(0.0,0.1) node {$v_{k+1}$};
\draw[above] (v1)++(0.0,0.1) node {$v_{1}$};
\draw[above] (w)++(0.0,0.1) node {$w$};
\draw (1,2) node {$T_0^*$};
\draw (5,2.5) node {$T_0$};
\draw [decorate, decoration = {brace}] (6.2,1.3) --  (1.8,1.3);
\draw (4,0.9) node {$P$};
\draw (7.5,0) node {$T$};
\end{scope}

\begin{scope}[yshift = -3.5cm, xshift = 0cm]
\path (6,1.5) coordinate (v1);
\path (4,1.5) coordinate (vk-1);
\path (3,1.5) coordinate (vk);
\path (2,1.5) coordinate (vk+1);
\draw (vk+1) -- (vk) -- (vk-1);
\draw (vk-1) -- (4.5,1.5);
\draw (5.5,1.5) -- (v1);
\draw (v1)  [fill=white] circle (\vr);
\draw (vk-1)  [fill=white] circle (\vr);
\draw (vk)  [fill=white] circle (\vr);
\draw (vk+1)  [fill=white] circle (\vr);
\draw[rounded corners=15pt] (0,0) rectangle ++(6.5,3);
\draw[rounded corners=15pt] (0.3,0.3) rectangle ++(2.2,2.4);
\draw[below] (v1)++(0.0,-0.1) node {$v_{1}'$};
\draw (5,2.5) node {$T_0'$};
\end{scope}

\end{tikzpicture}
\end{center}
\caption{The construction of the tree $T$ in Theorem~\ref{thm:double-c}}
\label{fig:construction}
\end{figure}

\begin{theorem}\label{thm:double-c}
Let $T_0$ be a tree of order $n\geq 7$ containing a proper pendant path $P$ of length $k$, where $v_{k+1}$ is its root. Let $T_0^*$ be the subtree of $T_0$ obtained by removing all the vertices of $P$ but $v_{k+1}$,
and let $T_0'$ be a copy of $T_0$ with the vertex $v_1'\in V(T_0')$ corresponding to $v_1$. Let $T$ be the tree obtained by joining the vertices $v_1$ and $v_1'$, and by attaching a new leaf $w$ at $v_1$.  See Fig.~\ref{fig:construction}. If $T_0^*$ is a partially transmission irrgular DBTM subtree of $T_0$ and $2n\in \Big(j^2+1,(j+1)^2\Big)$ with $j\in [k]$, then $T$ is  transmission irregular.
\end{theorem}

\proof
From the structure of $T$, we have $n(T)=2n+1$. Set $\Tr_T(v_1)=x$. Let $P=v_{k+1}v_k\cdots v_2v_1$ and let $P'=v_{k+1}'v_k'\cdots v_2'v_1'$ be the corresponding pendant path in $T_0'$. Then $\Tr_T(w)=x+2n-1$ and $\Tr_T(v_1')=x+1$ by Lemma~\ref{equal}.

Note that $P$ and $v_1v_1'P'$ are two internal paths of lengths $k$ and $k+1$, respectively, in $T$. Then $n_T(v_i)-n_T(v_{i+1})=n_T(v_i')-n_T(v_{i+1}')=2i+1$ for $i\in [k]$ from the structure of $T$. In view of Lemma \ref{equal} and $\Tr_T(v_1)=x$, we have $\Tr_T(v_i')=x+i^2$ and $\Tr_T(v_i)=x+i^2-1$ for $i\in [k+1]$, that is, $\Tr_T(P)=\{i^2-1:i\in [k+1]\}+x$ with $\Tr_T(P')=\Tr_T(V(P))+1$.

Next we consider the transmissions of vertices from $V(T)\setminus (V(P)\cup V(P')\cup\{w\})$. Set $V_0^*=V(T_0^*)$ and let $V_i$ be the set of vertices in $V_0^*$ at distance $i$ from $v_{k+1}$ in $T_0$. Then $V_0^* = \cup_{j=0}^{a} V_j$, where $a = \ecc_{T_0^*}(v_{k+1})$. For any edge $st\in E(T_0^*)$, without loss of generality, we may assume that $d_{T_0}(s,v_{k+1})>d_{T_0}(t,v_{k+1})$. Since $n_T(t)=n_{T_0}(t)+n+1$ and $n_T(s)=n_{T_0}(s)$, we have
\begin{equation}
n_T(t)-n_T(s)=n_{T_0}(t)-n_{T_0}(s)+n+1.
\label{eq1.1}
\end{equation}
Assume that $\Tr_{T_0}(u)-\Tr_{T_0}(v_{k+1})=h$ for any vertex $u\in V_j\subseteq V(T_0^*)$. By Lemma~\ref{equal} and~\eqref{eq1.1} we have $\Tr_T(u)=\Tr_{T}(v_{k+1})+h+j(n+1)$ with $\Tr_T(v_{k+1})=\Tr_{T_0}(v_{k+1})+(k+1)(n+1)+\Tr_{T_0}(v_1)$, that is, $\Tr_T(u)=\Tr_{T_0}(v_{k+1})+(k+j+1)(n+1)+\Tr_{T_0}(v_1)+h$. It follows that
\begin{equation}
\Tr_T(u)=\Tr_{T_0}(u)+j(n+1)+c
\label{eq1.2}
\end{equation}
for any vertex $u\in V_j$ where $c=(k+1)(n+1)+\Tr_{T_0}(v_1)$.

Note that $\Tr_{T_0}(T_0^*)$ is pairwise disjoint by the assumption. Let $s,t$ be arbitrary vertices of $T_0^*$.  Then $\Tr_T(s)\neq \Tr_T(t)$ for any $\{s,t\}\subseteq V_j$ with $j\in [a]$ because of~\eqref{eq1.2} and the fact that $\Tr_{T_0}(s)\neq \Tr_{T_0}(t)$. Assume that $s\in V_j$, $t\in V_\ell$ with $j,\ell\in [a]$ and $j\neq\ell$. Then  $\Tr_T(s)\neq \Tr_T(t)$ holds by~\eqref{eq1.2} and the fact that $T_0^*$ is a DBTM subtree of $T_0$. Therefore $\Tr_T(T_0^*)$ is pairwise disjoint.

Note that $V(T_0)=V(T_0^*)\cup V(P)$. Set $A_0=\Tr_T(T_0)$. From the structure of $T$, we have $\Tr_T(u)>\Tr_T(v_{k+1})=x+(k+1)^2-1$ for any vertex $u\in V(T_0^*)\setminus\{v_{k+1}\}$. Then $A_0$ is pairwise disjoint. By symmetry, we have $\Tr_T(T_0')=A_0+1$, which is also pairwise disjoint. In view of Lemma~\ref{equal} and the structure of $T$, the absolute value of the difference between any two elements in $A_0$ is greater than $1$. Therefore $A_0\cap A_1=\emptyset$ where $A_1=A_0+1$.  Recall that $\Tr_T(w)=x+2n-1$. Since $2n\in \Big(j^2+1,(j+1)^2\Big)$ with $j\in [k]$, we have $\Tr_T(w)\cap A=\emptyset$ with $A=A_0\cup A_1$.
\qed

Let $T$ be a transmission irregular chemical tree of order $n$ with a DBTM subtree $T_0$ obtained by removing all the non-root vertices of a pendant path of length $k$ such that $2n\in \Big(j^2+1,(j+1)^2\Big)$ with $j\in [k]$. By using the method in Theorem~\ref{thm:double-c}, we can construct another transmission irregular chemical tree of order $2n+1$.

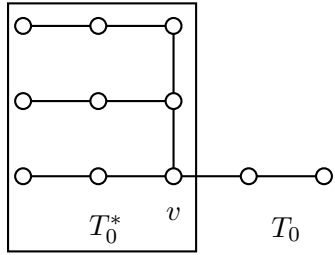
\begin{figure}[ht!]
\begin{center}
\begin{tikzpicture}[scale=1.0,style=thick]
\def\vr{3pt} 

\draw (-2,1) -- (-1,1) -- (0,1) -- (1,1) -- (2,1);
\draw (-2,2) -- (-1,2) -- (0,2) -- (0,1);
\draw (0,2) -- (0,3) -- (-1,3) -- (-2,3);

\draw (1,1)  [fill=white] circle (\vr);
\draw (0,1)  [fill=white] circle (\vr);
\draw (2,1)  [fill=white] circle (\vr);
\draw (-2,2)  [fill=white] circle (\vr);
\draw (0,2)  [fill=white] circle (\vr);
\draw (-1,1)  [fill=white] circle (\vr);
\draw (-1,2)  [fill=white] circle (\vr);

\draw (-0,3)  [fill=white] circle (\vr);
\draw (-1,3)  [fill=white] circle (\vr);
\draw (-2,3)  [fill=white] circle (\vr);

\draw (-2,1)  [fill=white] circle (\vr);

\draw (-2.2,0) rectangle (0.3,3.3);


\draw (1.5,0.3) node {$T_0$};

\draw (0,0.5) node {$v$};

\draw (-0.9,0.3) node {$T_0^*$};

\end{tikzpicture}
\end{center}
\caption{Tree $T_0$ with a non-DBTM  subtree $T_0^*$ rooted at $v$.}
\label{fig: ND}
\end{figure}

The condition that $T_0^*$ is a DBTM subtree of $T_0$ in Theorem~\ref{thm:double-c} is not necessary for obtaining a transmission irregular tree $T$. See an example in Fig.~\ref{fig: ND} of $T_0$ with a subtree $T_0^*$ rooted at vertex $v$ which is not DBTM. It is routine that the tree $T$, constructed from $T_0$ with the method in Theorem~\ref{thm:double-c}, is transmission irregular.

\section{Cycle-containing TI graphs}
\label{sec:cycle}

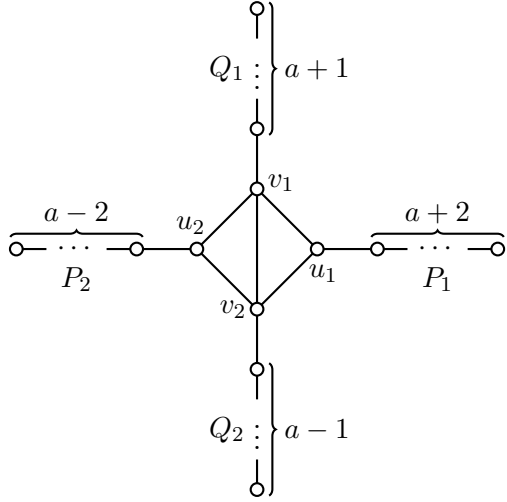
\begin{figure}[ht!]
\begin{center}
\begin{tikzpicture}[scale=0.8,style=thick]
\tikzstyle{every node}=[draw=none,fill=none]
\def\vr{3pt} 

\begin{scope}[yshift = 0cm, xshift = 0cm]
\path (1,0) coordinate (x1);
\path (2,0) coordinate (x2);
\path (4,0) coordinate (x3);
\path (-1,0) coordinate (x4);
\path (-2,0) coordinate (x5);
\path (-4,0) coordinate (x6);
\path (0,1) coordinate (x7);
\path (0,2) coordinate (x8);
\path (0,4) coordinate (x9);
\path (0,-1) coordinate (x10);
\path (0,-2) coordinate (x11);
\path (0,-4) coordinate (x12);
\draw (x1) -- (x7) -- (x4) -- (x10) -- (x1);
\draw (x7) -- (x10);
\draw (x1) -- (x2);
\draw (x4) -- (x5);
\draw (x7) -- (x8);
\draw (x10) -- (x11);
\draw (2,0) -- (2.5,0);
\draw (3.5,0) -- (4,0);
\draw (-2,0) -- (-2.5,0);
\draw (-3.5,0) -- (-4,0);
\draw (0,2) -- (0,2.5);
\draw (0,3.5) -- (0,4);
\draw (0,-2) -- (0,-2.5);
\draw (0,-3.5) -- (0,-4);
\draw (x1)  [fill=white] circle (\vr);
\draw (x2)  [fill=white] circle (\vr);
\draw (x3)  [fill=white] circle (\vr);
\draw (x4)  [fill=white] circle (\vr);
\draw (x5)  [fill=white] circle (\vr);
\draw (x6)  [fill=white] circle (\vr);
\draw (x7)  [fill=white] circle (\vr);
\draw (x8)  [fill=white] circle (\vr);
\draw (x9)  [fill=white] circle (\vr);
\draw (x10)  [fill=white] circle (\vr);
\draw (x11)  [fill=white] circle (\vr);
\draw (x12)  [fill=white] circle (\vr);
\draw (3,0) node {$\cdots$};
\draw (-3,0) node {$\cdots$};
\draw (0,3) node {$\vdots$};
\draw (0,-3) node {$\vdots$};
\draw (3.0,0.6) node {$a+2$};
\draw (-3.0,0.6) node {$a-2$};
\draw (1,3.0) node {$a+1$};
\draw (1,-3.0) node {$a-1$};
\draw (1.1,-0.4) node {$u_1$};
\draw (-1.1,0.4) node {$u_2$};
\draw (0.4,1.1) node {$v_1$};
\draw (-0.4,-1) node {$v_2$};
\draw (3,-0.5) node {$P_1$};
\draw (-3,-0.5) node {$P_2$};
\draw (-0.5,3) node {$Q_1$};
\draw (-0.5,-3) node {$Q_2$};
\draw [decorate, decoration = {brace}] (1.9,0.2) --  (4.1,0.2);
\draw [decorate, decoration = {brace}] (-4.1,0.2) -- (-1.9,0.2);
\draw [decorate, decoration = {brace}] (0.2,4.1) -- (0.2,1.9);
\draw [decorate, decoration = {brace}] (0.2,-1.9) --  (0.2,-4.1);
\end{scope}
\end{tikzpicture}
\end{center}
\caption{The graph $Z_0(a-1,a+1;a-2,a+2)$}
\label{fig:Z0}
\end{figure}

 Let $Z_0$ be the graph obtained from $K_4$ be removing one of its edges. For an integer $a\geq 2$, we denote by $Z_0(a-1,a+1;a-2,a+2)$ the graph obtained  from $Z_0$ by attaching a pendant path of length $a-1$ to a vertex of degree $3$, a pendant path of length $a+1$ at the other vertex of degree $3$, a pendant path of length $a-2$ at a vertex of degree $2$, and a pendant path of length $a+2$ at the other vertex of degree $2$, see Fig.~\ref{fig:Z0}.

\begin{theorem}\label{ZTI}
 If $a\geq 2$, then $Z_0(a-1,a+1;a-2,a+2)$ is TI if and only if $a$ is odd with $a\neq 1 \bmod 3$.
\end{theorem}

\proof
 Set $Z=Z_0(a-1,a+1;a-2,a+2)$. It is straightforward to check that $Z$ is not TI for $a=2$. Assume in the rest that $a\geq 3$.  Let $u_1$ and $u_2$ be the vertices in $Z$ of degree $3$, and let $P_1$ and $P_2$ be the paths of lengths $a+2$ and $a-2$ attached at $u_1$ and $u_2$, respectively. Let $v_1$ and $v_2$ be the vertices of degree $4$, and let $Q_1$ and $Q_2$ be the paths of lengths $a+1$ and $a-1$ attached to $v_1$ and $v_2$, respectively. See Fig.~\ref{fig:Z0} again.

Note that $n(Z)=4a+4$. Let $\Tr_Z(v_1)=x$. From the structure of $Z$ and Lemma~\ref{equal}, we have $\Tr_Z(v_2)=x+2$, $\Tr_Z(u_1)=a-2+x$, and $\Tr_Z(u_2)=a+6+x$. By Lemma~\ref{pend}, the transmissions of vertices lying on the pendant paths $P_1$, $P_2$, $Q_1$, and $Q_2$ including their roots, respectively  form the following sets:
\begin{align*}
A_1 & = \{(2j+1)a+(j-2)(j+1):\ j\in[a+2]_0\}+x, \\
A_2 & = \{(2j+1)a+(j+6)(j+1):\ j\in[a-2]_0\}+x, \\
B_1 & = \{2ja+j(j+1):j\in[a+1]_0\}+x, \\
B_2 & = \{2ja+j(j+5)+2:\ j\in[a-1]_0\}+x.
\end{align*}
Set $A=A_1\cup A_2$ and $B=B_1\cup B_2$. Therefore $Z$ is transmission irregular if and only if $A\cup B$ is pairwise disjoint.

If $a$ is even,  we select $(2k+1)a+(k-2)(k+1)\in A_1$ and $2ka+k(k+1)\in B_1$ with $k\in [a+1]$. Then $(2k+1)a+(k-2)(k+1)=2ka+k(k+1)$ if $k=\frac{a-2}{2}\in [a+1]$. It follows that $A\cap B\neq \emptyset$, that is, $Z$ is not transmission irregular. Next we turn to the case when $a$ is odd. Note that $A$ consists of odd numbers and $B$ consists of even numbers in this case. Therefore $A\cap B=\emptyset$ holds. To characterizing the TI property of $Z$ for odd $a$, it suffices to determine the condition of $a$ such that $A_1\cap A_2=\emptyset$ and $B_1\cap B_2=\emptyset$.

For any two elements $2sa+s(s+1)\in B_1$ with $s\in [a+1]_0$ and $2ta+t(t+5)+2\in B_2$ with $t\in [a-1]_0$, if $2sa+s(s+1)=2ta+t(t+5)+2$, then
\begin{equation}\label{eq0.0}
(s-t)(2a+s+t+1)=4t+2
\end{equation}
with $s-t>0$.  If $s-t=1$, then $2a+2t+2=4t+2$, that is, $t=a$. This contradicts the range of $t$. If $s-t=2$, we have $2(2a+2t+3)=4t+2$, which implies that $a=-1$. This is impossible. While $s-t\geq 3$, we have $s+t\geq 3$. From~\eqref{eq0.0} we have $(s-t)(2a+s+t+1)\geq 3(2a+4)=6a+12>4a-2=4(a-1)+2\geq 4t+2$. A clear contradiction occurs again. Therefore $B_1\cap B_2=\emptyset$ holds for any odd number $a$.

For any two elements $(2s+1)a+(s-2)(s+1)\in A_1$ with $s\in [a+2]_0$ and $(2t+1)a+(t+1)(t+6)\in A_2$ with $t\in [a-2]_0$, if $(2s+1)a+(s-2)(s+1)=(2t+1)a+(t+1)(t+6)$, then
\begin{equation}\label{eq0.1}
(s-t)(2a+t+s-1)=8t+8
\end{equation}
with $s>t$. If $s-t\geq 4$, then $s+t\geq 4$, which implies that $$(s-t)(2a+t+s-1)\geq 4(2a+3)=8a+12>8a-8\geq 8(a-2)+8\geq 8t+8.$$ Therefore, \eqref{eq0.1} does not hold. If $s-t=3$, then by \eqref{eq0.1} we have $3(2a+2t+2)=8t+8$, that is, $t=3a-1>a-2$, contradicting the fact that $t\in [a-2]_0$. If $s-t=2$, from~\eqref{eq0.1}, we have $2(2a+2t+1)=8t+8$, that is, $t=a-\frac{3}{2}$. This is impossible since $t$ is an integer. For $s-t=1$, similarly as above, we have $2a+2t=8t+8$, that is, $t=\frac{a-4}{3}$. Therefore, $A_1\cap A_2=\emptyset$ if and only if $a\neq 1 \bmod 3$. This completes the proof.
\qed

Denote by $K_4(k_1,k_2,k_3,k_4)$ the graph obtained  from the complete graph $K_4$ by  respectively attaching pendant paths of lengths $k_i\geq 0$, $i\in [4]$, to its vertices.

\begin{theorem}\label{k-4}
 If $a\geq 2$, then $K_4(a-2,a-1,a+1,a+2)$ is TI if and only if $a\neq 2 \bmod 3$.
\end{theorem}

\proof  Set $K=K_4(a-2,a-1,a+1,a+2)$.
For $a=2$, it can be easily checked that the vertex of degree $3$ has the same transmission as the vertex of degree $2$ adjacent to the vertex of degree $4$ at which a pendant path of length $4$ is attached. Therefore $K$ is not TI for $a=2$.

In the following assume that $a\geq 3$.   Let $v_1$, $v_2$, $v_3$, and $v_4$ be the vertices of degree $4$, and let $P^{(1)}$, $P^{(2)}$, $P^{(3)}$, and $P^{(4)}$ be respective attached paths of lengths $a+2$, $a+1$, $a-1$, and $a-2$. Let $\Tr(v_1)=x$. By Lemma~\ref{equal}, we have $\Tr(v_2)=x+1$, $\Tr(v_3)=x+3$, and $\Tr(v_4)=x+4$.  Note that  $n(K)=4a+4$. By Lemma~\ref{pend}, the set of transmissions of vertices on $P^{(1)}$ including $v_1$ is  $\{2ja+j(j-1):\ j\in[a+2]_0\}+x$. Similarly, the sets of transmissions of vertices on $P^{(2)}$, $P^{(3)}$, and $P^{(4)}$, respectively including $v_2$, $v_3$, and $v_4$, are  $\{2ja+j(j+1)+1:\ j\in[a+1]_0\}+x$, $\{2ja+(j+1)(j+4)-1:\ j\in[a-1]_0\}+x$ and $\{2ja+j(j+7)+4:\ j\in[a-2]_0\}+x$. For convenience, we set $A_1=\{2ja+j(j-1):\ j\in[a+2]_0\}$, $A_2=\{2ja+j(j+1)+1:\ j\in[a+1]_0\}$, $A_3=\{2ja+(j+1)(j+4)-1:\ j\in[a-1]_0\}$, $A_4=\{2ja+j(j+7)+4:\ j\in[a-2]_0\}$, and $A=\cup_{i=1}^{4}A_i$. Then $K$ is transmission irregular if and only if  the sets $A_i$, $i\in [4]$, are pairwise disjoint.

Note that each of $A_1$  and $A_4$ consists of even numbers and each of $A_2$ and $A_3$ consists of odd numbers. Clearly $A_p\cap A_q=\emptyset$ for any $p\in \{1,4\}$ and $q\in \{2,3\}$. Next we show that $A_2\cap A_3=\emptyset$. Otherwise, there are two elements $s=2ka+k(k+1)+1\in A_2$ and $t=2ja+(j+1)(j+4)-1\in A_3$ with $k\in [a+1]_0$, $j\in [a-1]_0$ and $s=t$. Clearly, we have $k>j$. Then $(k-j)(2a+k+j)=5j-k+2$, that is,
\begin{equation}\label{eq1}
  (k-j)(2a+k+j-1)=4j+2.
\end{equation}
If $k-j\geq 2$, then $(k-j)(2a+k+j-1)>4a+2$ and $4j+2\leq 4a-2$ since $j\in [a-1]_0$. Therefore $\eqref{eq1}$ does not hold. If $k-j=1$, we have $2a+2j=4j+2$ from~$\eqref{eq1}$. Then $a=j+1$, which is impossible since $j\in [a-1]_0$. Therefore $A_2\cap A_3=\emptyset$ follows immediately.

Now we determine the non-empty property of $A_1\cap A_4$. Choosing any two elements $s=2ka+k(k-1)\in A_1$ and $t=2ja+j(j+7)+4\in A_4$ with $k\in [a+2]_0$ and $j\in [a-2]_0$, we have $s-t=(k-j)(2a+k+j+7)-(8k+4)$. If $s=t$, then
\begin{equation}\label{eq2}
  (k-j)(2a+k+j+7)=8k+4
\end{equation}
with $k>j$. If $k-j\geq 4$,  then $(k-j)(2a+k+j+7)> 8a+28$ and $8k+4\leq 8a+20$ since $k\in [a+2]_0$. So $(\ref{eq2})$ does not hold. If $k-j=3$, then
$6a+6j+30=8j+28$ from~\eqref{eq2}, that is, $6a=2j-2$. Since $j\in [a-2]_0$, we have $6a\leq 2a-6$, contradicting the assumption $a\geq 3$. If $k-j=2$, similarly as above, we have $4a\leq 4a-6$ as a contradiction, again. If $k-j=1$, from Equality $(\ref{eq2})$, we have  $2a+2j+8=8j+12$, that is, $a=3j+2$. Therefore $s\neq t$, that is, $A_1\cap A_4=\emptyset$ if and only if $a\neq 2 \bmod 3$.
\qed

$Z_0(a-1,a+1;a-2,a+2)$ can be changed into $K_4(a-2,a-1,a+1,a+2)$  by adding an edge between the two vertices of degree $3$. By Theorems~\ref{ZTI} and~\ref{k-4}, if $a$ is an odd multiple of $3$, the TI property remains from $Z_0(a-1,a+1;a-2,a+2)$ to $K_4(a-2,a-1,a+1,a+2)$ by inserting a new edge.  In our last result we provide two sufficient conditions which guarantee that the transmission irregularity is preserved after inserting a new edge.

\begin{theorem}\label{1-cycle}
Let $G$ be a TI graph with $\Tr_G(v_1)>\Tr_G(v_2)>\Tr_G(v_3)$ as the first three largest transmissions. \begin{itemize}
\item [$(i)$]  If $v_1$, $v_2$, and $v_3$ lie on a pendant path $v_4v_3v_2v_1$ with natural adjacency relation, where $v_4$ is the root and $v_1$  is a  pendant vertex, then $G+v_2v_4$ is transmission irregular.
\item[$(ii)$]  If $v_1$ and $v_2$ are both pendant vertices with $v_1v_3\in E(G)$, $v_2$ and $v_3$ have a common neighbor and $\Tr_G(v_3)-1>\Tr_G(z)$  for any $z\in V(G)\setminus \{v_1, v_2, v_3\}$, then $G+v_2v_3$ is transmission irregular.
\end{itemize}
\end{theorem}

\proof
We first deal with (i).  For convenience, we set $G'=G+v_2v_4$ and  let $G_0$ be the subgraph of $G$ induced by $V_0$, where $V_0=V(G)\setminus\{v_1,v_2,v_3\}$.  Note that the vertices $v_2,v_3,v_4$ form a triangle in $G'$. For any vertex $u\in V_0$, we have $d_{G'}(u,v_i)=d_G(u,v_i)-1$ for $i\in [2]$, and $d_{G'}(u,w)=d_G(u,w)$ for any vertex $w$ in $(V_0\cup \{v_3\})\setminus\{v_1,v_2,u\}$. Therefore we have $\Tr_{G'}(u)=\Tr_G(u)-2$ for any vertex $u\in V_0$, that is, $\Tr_{G'}(G_0)=\Tr_G(G_0)-2$. Set $\Tr_G(v_4)=x$. Then, by Lemma~\ref{equal}, we have $\Tr_G(v_3)=x+n-6$, $\Tr_G(v_2)=x+2n-10$, and $\Tr_G(v_1)=x+3n-12$. Thus we have $\Tr(G)=\Tr_G(G_0)\cup\Big(\{n-6,2n-10,3n-12\}+x\Big)$.

From the structure of $G'$ and the above argument, we have $\Tr_{G'}(v_4)=x-2$, $\Tr_{G'}(v_3)=x+n-6$, $\Tr_{G'}(v_2)=x+n-7$, and $\Tr_{G'}(v_1)=x+2n-9$, which imply that
$$\Tr(G')=\Big(\Tr_G(G_0)-2\Big)\cup\Big(\{n-6,n-7,2n-9\}+x\Big).$$
From the assumption, we have $x+n-6>y$ for any $y\in \Tr_G(G_0)$, that is,
$$x+n-7>y-1>y-2$$
for any $y-2\in \Tr_G(G_0)-2$. Moreover, $\Big(\Tr_G(G_0)-2\Big)\cap\Big(\{n-6,n-7,2n-9\}+x\Big)=\emptyset$ with $\Tr_G(G_0)-2$ being pairwise disjoint. So $G'$ is transmission irregular as desired.

Now we turn to (ii). Assume that $v_4$ is the unique common neighbor of $v_2$ and $v_3$. Let $G^*=G+v_2v_3$ and $\Tr_G(v_4)=y$. By Lemma \ref{equal}, we have $\Tr_G(v_3)=y+n-4$, $\Tr_G(v_2)=y+n-2$, and $\Tr_G(v_1)=y+2n-6$. Note that $v_2$, $v_3$, and $v_4$ form a triangle in $G^*$. From the structure of $G^*$, we have  $d_{G^*}(u,w)=d_G(u,w)$ for any two vertices $u,w\in V(G)\setminus \{v_1,v_2,v_3\}$ and $d_{G^*}(u,z)=d_G(u,z)$ for any $z\in\{v_1,v_2,v_3\}$. Thus $\Tr_{G^*}(u)=\Tr_G(u)$ for any vertex $u\in V(G)\setminus \{v_1,v_2,v_3\}$ with $\Tr_{G^*}(v_3)=y+n-5$, $\Tr_{G^*}(v_2)=y+n-4$, and $\Tr_{G^*}(v_1)=y+2n-7$.

 Let $G_1 = G-\{v_1,v_2,v_3\}$. Then $\Tr_{G^*}(G)=\Tr_{G}(G_1)\cup\Big(\{n-5,n-4,2n-7\}+y\Big)$ from the above argument.  From the assumptions, $\Tr_G(G_1)$ is pairwise disjoint and $y+n-5>t$ for any $t\in \Tr_G(G_1)$. Therefore $G^*$ is transmission irregular.
\qed

\begin{figure}[ht!]
\begin{center}
\begin{tikzpicture}[scale=1.0,style=thick]
\def\vr{3pt} 

\draw (-1,1) -- (0,1) -- (1,1) -- (2,1) -- (3,1) -- (4,1) -- (5,1) -- (6,1) -- (7,1) -- (8,1);
\draw (-1,2) -- (0,2) -- (1,2) -- (2,2);
\draw (0,1) -- (0,0) -- (-1,0);
\draw (0,2) -- (0,1);
\draw (6,2) -- (6,1);
\draw (0,-2) -- (-1,-1) -- (0,-1) -- (0,0);
\draw (0,-1) -- (0,-2) -- (1,-2);

\draw (1,1)  [fill=white] circle (\vr);
\draw (0,1)  [fill=white] circle (\vr);
\draw (0,0)  [fill=white] circle (\vr);
\draw (-1,0)  [fill=white] circle (\vr);
\draw (2,1)  [fill=white] circle (\vr);
\draw (3,1)  [fill=white] circle (\vr);
\draw (4,1)  [fill=white] circle (\vr);
\draw (5,1)  [fill=white] circle (\vr);
\draw (6,1)  [fill=white] circle (\vr);
\draw (7,1)  [fill=white] circle (\vr);
\draw (8,1)  [fill=white] circle (\vr);
\draw (1,2)  [fill=white] circle (\vr);
\draw (6,2)  [fill=white] circle (\vr);
\draw (0,2)  [fill=white] circle (\vr);
\draw (-1,1)  [fill=white] circle (\vr);
\draw (-1,2)  [fill=white] circle (\vr);

\draw (-1,-1)  [fill=white] circle (\vr);
\draw (0,-1)  [fill=white] circle (\vr);
\draw (0,-2)  [fill=white] circle (\vr);
\draw (1,-2)  [fill=white] circle (\vr);

\draw (2,2)  [fill=white] circle (\vr);

\draw (-1,0.6) node {$19$};
\draw (1,0.6) node {$3$};
\draw (2,0.6) node {$8$};
\draw (3,0.6) node {$15$};
\draw (4,0.6) node {$24$};
\draw (5,0.6) node {$35$};
\draw (6,0.6) node {$48$};
\draw (6.4,1.95) node {$67$};
\draw (7,0.6) node {$65$};
\draw (8,0.6) node {$84$};

\draw (0,2.35) node {$13$};
\draw (-1,2.35) node {$32$};
\draw (1,2.35) node {$30$};
\draw (2,2.35) node {$49$};

\draw (0.3,0) node {$9$};
\draw (0.35,-1) node {$22$};
\draw (0.95,-1.7) node {$57$};
\draw (0.3,-1.75) node {$38$};
\draw (-1.4,0) node {$28$};
\draw (-1.4,-1) node {$39$};

\draw (0.2,0.6) node {$v$};

 \draw (0.2,1.25) node {$x$};

\draw (1.9,0) node {$G$};

\end{tikzpicture}
\end{center}
\caption{Graph $G$ satisfying the condition $(ii)$ in Theorem \ref{1-cycle}.}
\label{fig: G2}
\end{figure}
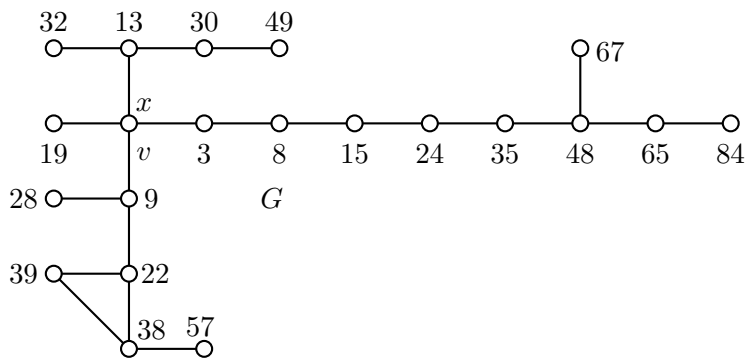

Two examples of graphs of order $21$ satisfying the conditions of $(i)$ and $(ii)$, respectively, in Theorem~\ref{1-cycle} are provided in Figs.~\ref{fig: G1} and~\ref{fig: G2} where a specific vertex $v$ is given with $\Tr_G(v)=x$ and $\Tr(G)=\{a_u:u\in V(G)\}+x$ for all the values of $a_u$ being labelled.  That is, the transmission of $v$ is $x$, and the transmission of every other vertex is equal to the sum of $x$ and the value next to the vertex.

\begin{figure}[ht!]
\begin{center}
\begin{tikzpicture}[scale=1.0,style=thick]
\def\vr{3pt} 

\draw (-3,1) -- (-2,1) -- (-1,1) -- (0,1) -- (1,1) -- (2,1) -- (3,1) -- (4,1) -- (4,2) -- (4,3) -- (5,3);
\draw (-2,1) -- (-2,0) -- (-1,0) -- (-1,1);
\draw (-1,3) -- (-1,2) -- (0,2) -- (0,1);
\draw (-1,-1) -- (0,-1) -- (0,0) -- (0,1) -- (1,2);
\draw (0,2) -- (-1,3) -- (-2,3);

\draw (1,1)  [fill=white] circle (\vr);
\draw (0,1)  [fill=white] circle (\vr);
\draw (0,0)  [fill=white] circle (\vr);
\draw (-1,0)  [fill=white] circle (\vr);
\draw (-2,0)  [fill=white] circle (\vr);
\draw (0,-1)  [fill=white] circle (\vr);
\draw (-1,-1)  [fill=white] circle (\vr);
\draw (-3,1)  [fill=white] circle (\vr);
\draw (2,1)  [fill=white] circle (\vr);
\draw (3,1)  [fill=white] circle (\vr);
\draw (4,1)  [fill=white] circle (\vr);
\draw (4,2)  [fill=white] circle (\vr);
\draw (4,3)  [fill=white] circle (\vr);
\draw (5,3)  [fill=white] circle (\vr);

\draw (1,2)  [fill=white] circle (\vr);
\draw (0,2)  [fill=white] circle (\vr);
\draw (-1,1)  [fill=white] circle (\vr);
\draw (-1,2)  [fill=white] circle (\vr);

\draw (-1,3)  [fill=white] circle (\vr);
\draw (-2,3)  [fill=white] circle (\vr);

\draw (-2,1)  [fill=white] circle (\vr);

 \draw (0.2,0.6) node {$x$};
\draw (1,0.6) node {$7$};
\draw (2,0.6) node {$16$};
\draw (3,0.6) node {$27$};
\draw (4,0.6) node {$40$};
\draw (4.4,1.9) node {$55$};
\draw (4.3,2.6) node {$72$};
\draw (5.1,2.6) node {$91$};
\draw (1.35,1.95) node {$19$};
\draw (0.35,1.9) node {$13$};
\draw (-0.6,3.1) node {$29$};
\draw (-2.4,3.1) node {$48$};
\draw (-1.4,2) node {$30$};

\draw (0.35,0) node {$15$};
\draw (0.35,-1) node {$32$};
\draw (-0.9,-1.35) node {$51$};

\draw (-1.1,1.3) node {$11$};
\draw (-2.1,1.3) node {$26$};
\draw (-1.1,-0.3) node {$28$};
\draw (-2.1,-0.3) node {$43$};
\draw (-3.1,1.3) node {$45$};

\draw (-0.2,1.2) node {$v$};

\draw (1.3,-0.6) node {$G$};

\end{tikzpicture}
\end{center}
\caption{Graph $G$ satisfying the condition $(i)$ in Theorem \ref{1-cycle}.}
\label{fig: G1}
\end{figure}
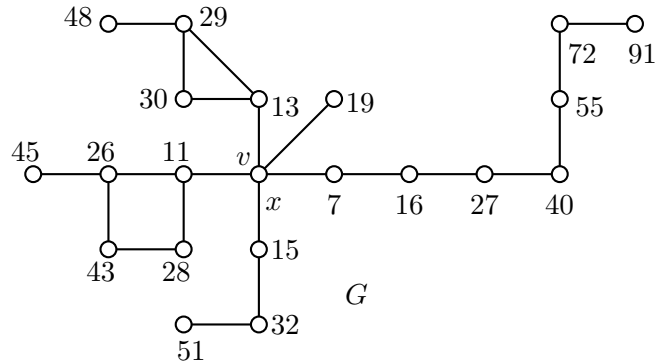

\section{Concluding remarks}
\label{sec:remark}

In this paper we prove the TI property of some chemical graphs and provide the method of construct new TI graphs.

Note that there are some transmission irregular chemical graphs of even order $n=4a+4$ with $a\neq 2 \bmod 3$.  Combining this fact with Theorem~\ref{C-TI}, we pose the following problem.

\begin{problem}
 Does there exist a TI chemical graph of every even order?
\end{problem}

 Note that Theorem \ref{C-TI} implicitly provides a method constructing  TI chemical trees from known TI ones by attaching a pendant vertex at each of their leaves. Theorem~\ref{1-cycle} naturally leads to the following two problems.

\begin{problem}
 Establish additional methods for constructing TI graphs from known TI graphs.
\end{problem}

\begin{problem}\label{p2}
 Characterize TI chemical graphs $G$ which preserve TI property after joining two nonadjacent vertices.
\end{problem}

\section*{Acknowledgments}

 K.~Xu and J.~Tian were partially supported by NNSF of China (Grant No.\ 12271251). S.~Klav\v{z}ar was partially supported by the Slovenian Research Agency (ARRS) under the grants P1-0297, J1-2452, and N1-0285.

\end{document}